\documentclass[11pt,letterpaper]{amsart}
\usepackage{geometry}
\usepackage{hyperref,mathrsfs,url,setspace,mathpazo}
\newtheorem{theorem}{Theorem}
\newtheorem{proposition}{Proposition}
\newtheorem{lemma}{Lemma}
\newtheorem{corollary}{Corollary}
\theoremstyle{remark}
\newtheorem{remark}{Remark}
\newcommand{\C}{\mathbb{C}}
\newcommand{\D}{\Omega}
\newcommand{\ep}{\varepsilon}
\newcommand{\Dc}{\overline{\Omega}}
\newcommand{\dbar}{\overline{\partial}}
\onehalfspace

\title{Axler-Zheng type theorem on a class of domains in $\C^n$}
\author{\u{Z}eljko \u{C}u\u{c}kovi\'c}
\author{S\"{o}nmez \c{S}ahuto\u{g}lu}
\email{Zeljko.Cuckovic@utoledo.edu, Sonmez.Sahutoglu@utoledo.edu}
\address{University of Toledo, Department of Mathematics \& Statistics, 
Toledo, OH 43606, USA}


\date{\today}
\begin{document}

\begin{abstract}
We prove a version of Axler-Zheng's Theorem  on smooth bounded pseudoconvex domains in 
$\C^n$ on which the $\dbar$-Neumann operator is compact. 
\end{abstract}
\maketitle

When studying operators on spaces of holomorphic spaces, one of the important
questions is about the compactness of operators. It has been known that on
Bergman spaces, it is the Berezin transform that plays the crucial role in
characterizing compact operators. 
Let $\widetilde{S}$ denote the Berezin transform of $S$ (defined below). 
Axler-Zheng Theorem \cite{AxlerZheng98} states that if $S$ is a finite sum of
finite products of Toeplitz operators on $A^2(\mathbb{D})$ with symbols in
$L^{\infty}(\mathbb{D})$, then $S$ is compact if and only if
$\widetilde{S}(z)\to 0$ as $|z|\to 1^-$.
This result has been completed by Su\'arez
\cite{Suarez07} to include operators in the Toeplitz algebra even on the unit
ball of $\mathbb{C}^n$. Engli\v{s} \cite{Englis99} prove the analog of
Axler-Zheng Theorem for Fock (the Segal-Bargmann) space and irreducible
bounded symmetric domains in $\C^n$. We note that Axler-Zheng Theorem is 
also true on the polydisk as the irreducibilty condition in \cite{Englis99} can be 
removed (see  \cite[pg. 232]{ChoeKooLee09}). Raimondo \cite{Raimondo00} 
found a proof of Axler-Zheng Theorem for different domains in $\mathbb{C}$. 
Recently, Mitkovski and Wick extended Axler-Zheng Theorem to the Bergman 
spaces $A^p(\mathbb{D}^n)$ and with Su\'arez  they treated weighted Bergman 
spaces on the unit ball in $\C^n$ (see
\cite{MitkovskiSuarezWick13,MitkovskiWick}). 

The purpose of this paper is to provide a version of Axler-Zheng Theorem for
quite general domains in $\C^n$, namely for smooth bounded pseudoconvex domains
in $\C^n$. This setting provides a framework for the use of the methods from 
several complex variables, in particular the $\dbar$-Neumann problem. However,
that forces symbols to be continuous up to the boundary. Before, we state our
main result, we will introduce notation and basic definitions. 

Let $\D$ be a smooth bounded pseudoconvex domain in $\C^n$. The Bergman space
$A^2(\D)$ is the closed subspace of $L^2(\D)$(with Lebesgue measure $dV$) 
consisting of holomorphic functions on $\D$. Let $P:L^2(\D)\to A^2(\D)$ denote
the Bergman projection. For $\phi\in L^{\infty}(\D)$, the Toeplitz operator
with symbol $\phi$ is defined as $T_{\phi}f=P(\phi f)$, for $f\in A^2(\D)$. Let 
$K$ denote the Bergman kernel for $\D.$ Then $K$ is the unique function on
$\D\times \D$ with the following property: $f(z)= \langle f,K(z,.)\rangle$, for
every $f\in A^2(\D)$ and $z\in \D$. Then
\[k_z(w)=\frac{K(w,z)}{\sqrt{K(z,z)}}\] 
is called the normalized Bergman kernel. One can check that $\|k_z\|=1.$ Let
$\mathscr{B}(A^2(\D))$ denote the algebra of bounded linear operators on
$A^2(\D).$ For $S\in \mathscr{B}(A^2(\D))$, the Berezin transform of $S$ is
defined as the function 
\[\widetilde{S}(z)=\langle Sk_z,k_z\rangle \text{ for } z\in \D.\] 
We define the Berezin transform of a function $\phi\in L^{\infty}(\D)$ as
$\widetilde{\phi}= \widetilde{T}_{\phi}.$ 

Let $\mathscr{T}(\Dc)$ denote the norm closed subalgebra of 
$\mathscr{B}(A^2(\D))$ generated by $\{T_{\phi}:\phi\in C(\Dc)\}.$  This algebra
has been studied by Coburn \cite{Coburn73} in the case of $\D$ is the open unit
disk in $\mathbb{C}$ and by Axler, Conway, and McDonald 
\cite{AxlerConwayMcDonald82} in case of $\D$ is a nonempty bounded domain in
$\C$.  Now we state our main result. 

\begin{theorem}\label{ThmMain}
Let $\D$ be a smooth bounded pseudoconvex domain in $\C^n$ on which 
$\dbar$-Neumann operator is compact. Assume that $T\in \mathscr{T}(\Dc)$
and $\lim_{z\to p}\widetilde{T}(z)=0$ for all $p\in b\D.$ Then $T$ is compact.
\end{theorem}

Our method uses standard arguments from the $\dbar$-Neumann problem. When the
domain $\D$ is bounded pseudoconvex, as a result of H\"{o}rmander's work, the
complex  Laplacian $\dbar^*\dbar+\dbar\dbar^*$ (here $\dbar^*$ is the Hilbert
space adjoint of $\dbar$)  has a bounded inverse $N$, known as the 
$\dbar$-Neumann operator,  on the square integrable $(0,1)$-forms on $\D$. We
refer the reader to  \cite{ChenShawBook,StraubeBook} for more information on the
$\dbar$-Neumann problem.

\begin{remark}
The proof of Theorem \ref{ThmMain} breaks if the operator $T$  is not
in $\mathscr{T}(\Dc)$.  In fact, if the domain $\D$ is the unit ball in $\C^n$, 
Theorem \ref{ThmMain} is true whenever $T$ is in the Toeplitz algebra. So in
that case $T$ does not have to be in  $\mathscr{T}(\Dc)$ (see \cite{Suarez07}).
It would be  interesting to know if the condition $T\in \mathscr{T}(\Dc)$ is 
necessary in our theorem.   
\end{remark}

\begin{remark}\label{Rmk1}
The class of domains where the $\dbar$-Neumann operator is compact is very
large; yet it excludes certain domains with analytic disks in the boundary. A
set $X$ is said to have an analytic disc if there exists a nonconstant
holomorphic mapping $f:\mathbb{D}\to X$. The $\dbar$-Neumann operator is
compact on strongly pseudoconvex domains and more generally on domains
with property $(P)$ (or $(\widetilde{P})$) (see \cite{Catlin84,McNeal02} and
\cite[Theorem 4.8]{StraubeBook}). For example, a bounded convex domain satisfies
Property $(P)$ if and only if there is no analytic disk in its boundary (see
\cite{FuStraube98} and \cite[Theorem 4.26]{StraubeBook}). Hence the next
corollary follows immediately. 
\end{remark}

\begin{corollary}\label{Cor1}
 Let $\D$ be a  smooth bounded convex domain in $\C^n$ whose boundary has
no analytic disk. Assume that $T\in \mathscr{T}(\Dc)$ and
$\lim_{z\to p}\widetilde{T}(z)=0$ for all $p\in b\D.$ Then $T$ is compact. 
\end{corollary}

\begin{remark}
As mentioned in Remark \ref{Rmk1}, Theorem \ref{ThmMain}  applies to a large
class of domains in $\C^n$. However, it excludes some simple domains  such as
the polydisk as the $\dbar$-Neumann operator is not compact on
the polydisk \cite{Salinas91} (see also \cite{Krantz88,FuStraube98}).
In case $\D$ is a smooth bounded strongly pseudoconvex in $\C^n$, Theorem
\ref{ThmMain} is a special case of \cite[Theorem 2.3]{ArazyEnglis01}. 
\end{remark}

\begin{corollary}\label{Cor2}
 Let $\D$ be a  smooth bounded convex domain in $\C^2$ and
$\phi_j\in C^1(\Dc)$ for $1\leq j\leq m$. Assume that $\phi_j\circ
f$ is holomorphic for all $j$ and any holomorphic function $f:\mathbb{D}\to
b\D$. Then the following are equivalent 
\begin{itemize}
 \item[i.] $T=T_{\phi_1}\cdots T_{\phi_m}$ is compact,
 \item[ii.] $\lim_{z\to p}\widetilde{T}(z)=0$ for all $p\in b\D,$
 \item[iii.] the product $\phi_1\phi_2\cdots \phi_m=0$ on $b\D$.
\end{itemize}
\end{corollary}
\begin{remark}
Corollary \ref{Cor2} applies to finite sums of finite products of Toeplitz
operators as well.  Let $\D$ be a smooth bounded convex domain in $\C^2$ and
$\phi_j\in C(\Dc)$ for $1\leq j\leq m$. Assume that $\phi_j\circ f$ is
holomorphic for all $j$ and any holomorphic function $f:\mathbb{D}\to b\D$ and
that $T$ is a finite sum of finite products of Toeplitz operators $T_{\phi_j}$. 
Then compactness of $T$ is equivalent to the continuous function $\phi$ being
zero on $b\D$ where $\phi$ is the corresponding finite sum of finite products of
the symbols $\phi_j$.
\end{remark}

The techniques in this paper also give us some results about the essential
norms of Toeplitz operators with symbols continuous up to the boundary. 
Let $T:X\to Y$ be a bounded linear operator between two Banach spaces
$X$ and $Y$.  Then the essential norm $\|T\|_e$ of $T$ is defined as
\[\|T\|_e=\inf \{ \|T-S\|: S :X\to Y \text{ is a compact linear operator}\}.\]
We remind the reader that for any $\psi\in C(\Dc)$ its sup norm on the
boundary is defined as 
\[\|\psi\|_{L^{\infty}(b\D)}=\sup\{|\psi(z)|:z\in b\D\}.\]

For  a bounded domain $\D$ in $\C$,  Axler, Conway, and McDonald
\cite{AxlerConwayMcDonald82} obtained the essential spectrum and the essential
norm of an arbitrary Toeplitz operator with a symbol in $C(\Dc)$. The following
result is an extension of their result about the essential norm.

\begin{proposition}\label{PropEssentialNorm}
 Let $\D$ be a bounded domain in $\C^n$ and $\psi\in C(\Dc).$
Then $\|T_{\psi}\|_e \leq \|\psi\|_{L^{\infty}(b\D)}.$ Furthermore, if $\D$ is
smooth bounded pseudoconvex in $\C^n$ such that the set of strongly pseudoconvex
points is dense in $b\D$ then $\|T_{\psi}\|_e = \|\psi\|_{L^{\infty}(b\D)}.$
\end{proposition}

We get the following immediate corollaries due to the following facts:  
compactness of the $\dbar$-Neumann operator implies that the set of strongly
pseudoconvex points is dense in the boundary in relative topology 
\cite[Corollary 1]{StraubeSahutoglu06}; and in case of convex domains, 
absence of analytic disks in the boundary is equivalent to compactness of the 
$\dbar$-Neumann operator \cite[Theorem 1.1]{FuStraube98}.
\begin{corollary}
Let $\D$ be smooth bounded pseudoconvex domain in $\C^n$ and $\psi\in C(\Dc)$
such that the $\dbar$-Neumann operator is compact. Then 
$\|T_{\psi}\|_e = \|\psi\|_{L^{\infty}(b\D)}.$
\end{corollary}

\begin{corollary}
Let $\D$ be smooth bounded convex in $\C^n$ with no analytic
disks in the boundary and $\psi\in C(\Dc)$. Then 
$\|T_{\psi}\|_e = \|\psi\|_{L^{\infty}(b\D)}.$
\end{corollary}

\section*{Proofs}
To prove Theorem \ref{ThmMain} we use the following lemma.

\begin{lemma} \label{LemDelta}
Let $\D$ be a smooth bounded pseudoconvex domain in 
$\C^n, \phi\in C(\Dc),$ and $p\in b\D$ be a strongly pseudoconvex point. Then
$\lim_{z\to p} \widetilde{\phi}(z)= \phi(p).$  
\end{lemma}
\begin{proof}
Without loss of generality we will assume that $\phi$ is real valued. Let
$B(p,r)$ denote the open ball centered at $p$ with radius $r.$ Continuity of
$\phi$ on $\Dc$ implies that for every $\ep>0$ there exists $\delta_1>0$ such
that $|\phi(w)-\phi(p)|<\ep$ for $w\in \Dc\cap B(p,\delta_1).$
Smoothness of the domain implies that there exists $0<r<\delta_1$ so
that every point in $B(p,r)\cap b\D$ is strongly pseudoconvex. 
Then $\Gamma_1=b\D\setminus \overline{B(p,r)}$ and 
$\Gamma_2=b\D\cap B(p,r/2)$ are disjoint open subsets of $b\D.$ Then  
\cite[Theorem 2]{Bell86} (see also \cite[Theorem 2]{Boas87}) implies that the
Bergman kernel $K$ extends smoothly to $\Gamma_1\times \Gamma_2$. Since
$K(z,z)\to \infty$ as $z\to p$ (\cite{Pflug75}, see also 
\cite[Theorem 6.1.7]{JarnickiPflugBook1}),  the normalized Bergman kernel
$|k_z|^2$ converges uniformly to 0 on $\Dc\setminus B(p,\delta_1)$ as $z\to p.$ 
 Then there exists $0<\delta_2<r/2$ such that  
$z\in \Dc\cap B(p,\delta_2)$ implies that $|k_z(w)|^2< \ep$ for all 
$w\in \Dc\setminus B(p,\delta_1).$ First we will estimate $\widetilde{\phi}(z)$
 from above for $z\in \Dc\cap B(p,\delta_2)$. Let
$M=\int_{\D}|\phi(w)|dV(w)+1.$ Then
\begin{align*}
\widetilde{\phi}(z)&=\int_{\D}\phi(w)|k_z(w)|^2dV(w)\\
&\leq \int_{\D\setminus B(p,\delta_1)} |\phi(w)| |k_z(w)|^2 dV(w)
+(\phi(p)+\ep)\int_{\D\cap B(p,\delta_1)} |k_z(w)|^2dV(w) \\
&\leq \varepsilon \int_{\D\setminus B(p,\delta_1)}|\phi(w)|dV(w)
+(\phi(p)+\ep)\int_{\D\cap B(p,\delta_1)} |k_z(w)|^2dV(w)\\
&\leq \varepsilon\int_{\D}|\phi(w)|dV(w)
+(\phi(p)+\ep)\left(1-\int_{\D\setminus B(p,\delta_1)}|k_z(w)|^2dV(w)\right)\\
&\leq \phi(p)+\ep M+|\phi(p)|\int_{\D\setminus B(p,\delta_1)}|k_z(w)|^2dV(w)\\
&\leq \phi(p)+\ep(M+|\phi(p)|V(\D)).
\end{align*}
In the last step we used that fact that $\int_{\D\setminus
B(p,\delta_1)}|k_z(w)|^2dV(w)\leq \varepsilon V(\D)$, where
$V(\D)$ denotes the volume of $\D.$
Hence for $z\in \Dc\cap B(p,\delta_2)$ we have 
\begin{align}\label{Eq1}
 \widetilde{\phi}(z)-\phi(p)  \leq \varepsilon (M+V(\D)|\phi(p)|).
\end{align}
Similarly we estimate $\widetilde{\phi}(z)$ from below as follows: 
\begin{align*}
\widetilde{\phi}(z)&=\int_{\D}\phi(w)|k_z(w)|^2dV(w)\\
&\geq -\int_{\D\setminus B(p,\delta_1)} |\phi(w)| |k_z(w)|^2 dV(w)
+(\phi(p)-\ep)\int_{\D\cap B(p,\delta_1)} |k_z(w)|^2dV(w) \\
&\geq -\varepsilon \int_{\D\setminus B(p,\delta_1)}|\phi(w)|dV(w)
+(\phi(p)-\ep)\int_{\D\cap B(p,\delta_1)} |k_z(w)|^2dV(w)\\
&\geq -\varepsilon\int_{\D}|\phi(w)|dV(w)
+(\phi(p)-\ep)\left(1-\int_{\D\setminus B(p,\delta_1)}|k_z(w)|^2dV(w)\right)\\
&\geq \phi(p)-\ep M-|\phi(p)|\int_{\D\setminus B(p,\delta_1)}|k_z(w)|^2dV(w)\\
&\geq \phi(p)-\ep(M+|\phi(p)|V(\D)).
\end{align*}
Hence for $z\in \Dc\cap B(p,\delta_2)$ we have 
\begin{align}\label{Eq2}
 \widetilde{\phi}(z)-\phi(p)  \geq -\varepsilon (M+V(\D)|\phi(p)|).
\end{align}
 Combining \eqref{Eq1} and \eqref{Eq2}  we showed that for $\ep>0$
there exists $\delta_2>0$ so that 
\begin{align*}
 -\varepsilon M \leq \widetilde{\phi}(z)-\phi(p) \leq \varepsilon M
\end{align*}
 for $z\in \Dc\cap B(p,\delta_2).$ Therefore, 
 $\lim_{z\to p}\widetilde{\phi}(z)=\phi(p).$
\end{proof}

\begin{remark}
The Berezin transform of a function that is continuous up to the boundary is
continuous  on strongly pseudoconvex points in the boundary.
\end{remark}

\begin{corollary}
 Let $\D$ be a  smooth bounded pseudoconvex domain in $\C^n$ such
that the set of strongly pseudoconvex points is dense in $b\D.$ Assume that
$\phi\in C(\Dc)$ and $\lim_{z\to p}\widetilde{\phi}(z)=0$ for all strongly
pseudoconvex points $p\in b\D$. Then $T_{\phi}$ is compact.
\end{corollary}

Next we prove Proposition \ref{PropEssentialNorm}. 
\begin{proof}[Proof of Proposition \ref{PropEssentialNorm}]
 Let  us fix $\alpha>\|\psi\|_{L^{\infty}(b\D)}\geq 0.$ We will show that
there exists a compact operator $S$ such that $\|T_{\psi}-S\|<\alpha.$ Since
$b\D$ is compact and $|\psi|<\alpha$ on $b\D$ there exists an open neighborhood
$U$ of $b\D$ such that $|\psi|< \alpha $ on $U\cap \D.$ Let us choose a smooth
cut-off function $\chi\in C^{\infty}_0(\D)$ such that $0\leq \chi\leq 1,\chi=1$
on a neighborhood of $\D\setminus U,$ and $\chi=0$ on $b\D.$  Then
$S=T_{\chi\psi}$ is a compact operator because $\chi\psi=0$ on $b\D.$
Furthermore, 
\[\|T_{\psi}\|_e\leq \|T_{\psi}-S\|\leq 
\sup\{|(1-\chi(z))\psi(z)|:z\in \Dc\}< \alpha.\]
Since $\alpha$ is arbitrary we get $\|T_{\psi}\|_e  \leq
\|\psi\|_{L^{\infty}(b\D)}.$

To prove the second part assume that $p\in b\D$ is a strongly pseudoconvex
point and $\|T_{\psi}\|_e=\alpha_0.$ Then for every $\ep>0$ there exists a
compact operator $K$ such that $\|T_{\psi}-K\| < \alpha_0+\ep.$ Applying the
Berezin transform we get $|\langle (T_{\psi}-K)k_z,k_z\rangle| < \alpha_0
+\ep$ for all $z\in \D.$ If we let $z$ converge to $p$, compactness of $K$
implies that $\langle Kk_z,k_z\rangle$ converges to zero and since $p$ is
strongly pseudoconvex Lemma \ref{LemDelta} implies  
\[|\psi(p)|=\lim_{z\to p}|\widetilde{\psi}(z)| <\alpha_0+\ep.\] 
Therefore,  we showed that $|\psi(z)| \leq \alpha_0$ for any strongly
pseudoconvex
point $z\in b\D$. Then continuity of $\psi$ on $b\D$ and the assumption
that the set of strongly pseudoconvex points is dense in the boundary imply
that $\|\psi\|_{L^{\infty}(b\D)}\leq \|T_{\psi}\|_e$. 
Hence the proof of Proposition \ref{PropEssentialNorm} is complete.
\end{proof}

 In the proof of Theorem \ref{ThmMain} we use Hankel operators. The Hankel operator 
 $H_{\phi}:A^2 (\D)\to L^2 (\D)$ with a symbol $\phi \in L^{\infty} (\D)$ is defined by 
 $H_{\phi}(f) = (I- P )(\phi f)$. Kohn's formula $P = I- \dbar^* N\dbar$ 
implies that $H_{\phi} (f ) = \dbar^* N\dbar (\phi f )$ for $f \in A^2 (\D)$ and
$\phi\in C  (\Dc)$. The following standard fact will be useful: if $N$ is a
compact operator on $\D$ then $H_{\phi}$ is compact for all $\phi\in C(\Dc)$ 
(see \cite[Propositions 4.1 and 4.2]{StraubeBook}). 

Now we are ready to prove the main theorem.

\begin{proof}[Proof of Theorem \ref{ThmMain}]
First we will prove the theorem in case of $T$ is a finite product of Toeplitz
operators with symbols continuous on $\Dc.$ Assume that
$T=T_{\phi_m}T_{\phi_{m-1}}\cdots T_{\phi_1}$ such that $\phi_j\in C(\Dc)$ for all $j$ and 
$\lim_{z\to p}\widetilde{T}(z)=0$ for all $p\in b\D.$ Then $\widetilde{T}$ has a
continuous extension on $\Dc$ and $\widetilde{T}=0$ on $b\D.$  

 Using the formula
$H_{\overline{\psi}_1}^*H_{\psi_2}=T_{\psi_1\psi_2}-T_{\psi_1}T_{\psi_2}$
inductively one can show that 
\begin{align}\nonumber 
T_{\phi_m}T_{\phi_{m-1}}\cdots T_{\phi_2}T_{\phi_1} 
= &T_{\phi_m\phi_{m-1}\cdots\phi_2\phi_1} - 
\Big( T_{\phi_m}T_{\phi_{m-1}} \cdots T_{\phi_3} 
H^*_{\overline{\phi}_2}H_{\phi_1}  \\ \label{EqnReduction} 
&+ T_{\phi_m}T_{\phi_{m-1}} \cdots T_{\phi_4}
H^*_{\overline{\phi}_3}H_{\phi_2\phi_1} 
+\cdots +H^*_{\overline{\phi}_m}H_{\phi_{m-1}\cdots\phi_2\phi_1} \Big) .
\end{align}
Compactness of the $\dbar$-Neumann operator implies that Hankel operators
with symbols continuous on the closure are compact
\cite[Proposition 4.1 and 4.2]{StraubeBook}.
Therefore, we have 
\[T_{\phi_m}T_{\phi_{m-1}}\cdots T_{\phi_1} 
=T_{\phi_m\phi_{m-1}\cdots\phi_1}+K\]
where   $K$ is a compact operator. 

Let $\phi=\phi_m\phi_{m-1}\cdots\phi_1$ on $\Dc.$ Then $T=T_{\phi}+K$  and 
\[ \widetilde{T}=\widetilde{T}_{\phi}+\widetilde{K}.\] 
Since $K$ is compact we have $\widetilde{K}=0$ on $b\D.$ Hence,
$\widetilde{T}_{\phi}=0$ on $b\D.$ However, Lemma \ref{LemDelta} 
implies that $\phi=\widetilde{\phi}=\widetilde{T}_{\phi}=0$ on strongly
pseudoconvex points. Since $\dbar$-Neumann operator is compact, \cite[Corollary
1]{StraubeSahutoglu06} implies that the set of strongly pseudoconvex points is dense in the
boundary (in the relative topology). This fact together with the continuity of $\phi$ implies
that $\phi=0$ on $b\D.$ Hence, $T_{\phi}$ and in turn $T$  are compact. 

In case $T$ is a finite sum of finite products of Toeplitz operators with
symbols continuous on $\Dc$ equation \eqref{EqnReduction} implies that there
exists $\phi\in C(\Dc)$ and a compact operator $K$ such that $T=T_{\phi}+K.$
Then by the arguments in the previous paragraph we conclude that $\phi=0$ on
$b\D.$ Hence, in this case too, $T$ is compact. 

Finally we will prove the Theorem \ref{ThmMain} for $T\in \mathscr{T}(\Dc).$
Using \eqref{EqnReduction} we can assume that for every $\varepsilon>0$ there
exist $\psi_{\varepsilon}\in C(\Dc)$ and a compact operator $K_{\varepsilon}$
such that $\|T-T_{\psi_{\varepsilon}}+K_{\varepsilon}\|<\varepsilon.$  Then for
$z\in \D$ we have 
\begin{align*}
 \left| \widetilde{T}(z)-\widetilde{T}_{\psi_{\varepsilon}}(z)
+\widetilde{K}_{\varepsilon}(z)\right| 
&= \left|\langle T(k_z)-T_{\psi_{\varepsilon}}(k_z)+K_{\varepsilon}(k_z),
k_z\rangle\right|\\
&\leq \|T-T_{\psi_{\varepsilon}}+K_{\varepsilon}\| \\
&< \varepsilon 
\end{align*}
Since we assumed that $\lim_{z\to p}\widetilde{T}(z)=0$ for all  $p\in b\D$ and
$K_{\varepsilon}$ is compact, we have
\[\limsup_{z\to
p}\left|\widetilde{T}_{\psi_{\varepsilon}}(p)\right|<\varepsilon\]
for every $p\in b\D.$ However, Lemma \ref{LemDelta} implies that 
$|\psi_{\varepsilon}(p)|=\lim_{z\to p}\left|\widetilde{\psi_{\varepsilon}}
(z)\right|<\varepsilon$ for any strongly pseudoconvex point $p\in b\D.$ 
Then using continuity of $\psi_{\varepsilon}$
on $\Dc$ together with the fact that strongly pseudoconvex points are dense in
the boundary, we deduce that $\left|\psi_{\varepsilon}\right|<\varepsilon$ on
$b\D.$ Proposition \ref{PropEssentialNorm} implies that there exists a
compact operator $S_{\varepsilon}$ such that 
$\|T_{\psi_{\varepsilon}}-S_{\varepsilon}\|<\varepsilon.$  Then 
\begin{align*}
\|T-S_{\varepsilon}+K_{\varepsilon}\|\leq
\|T-T_{\psi_{\varepsilon}}+K_{\varepsilon}\|
+\|T_{\psi_{\varepsilon}}-S_{\varepsilon} \| < 2\varepsilon.  
\end{align*}
Therefore, $T$ is in the operator norm closure of subspace of compact
operators in the operators. Since this subspace is closed in the
operator norm topology we conclude that $T$ is compact. 
\end{proof}

To prove Corollary \ref{Cor2} we will need the following result. We note that
even though the result below was stated for symbols that are smooth up to the
boundary, its proof goes through for symbols that are continuously 
differentiable on the closure of the domain.

 \begin{theorem}[{\cite[Theorem 3]{CuckovicSahutoglu09}}] \label{ThmCS09}
Let $\D$ be a smooth bounded convex domain in $\C^2$ and 
$\phi \in C^1 (\Dc)$. If $\phi\circ f$ is holomorphic for any
holomorphic $f :\mathbb{D} \to b\D$, then $H_{\phi}$ is compact.
 \end{theorem}

\begin{proof}[Proof of Corollary \ref{Cor2}]
The implication i. $\Rightarrow$ ii. is a result of compactness of $T$ and the
fact that $k_z$ converges to zero weakly as $z$ converges to $b\D$. 

To prove ii. $\Rightarrow$ iii.  we use Theorem \ref{ThmCS09} to conclude that $H_{\phi_j}$
is compact for every $j$. This fact together with \eqref{EqnReduction} implies
that $T=T_{\phi}+K$ where $\phi=\phi_1\phi_2\cdots\phi_n\in C(\Dc)$ and $K$ is a
compact operator. Then ii. implies that $\widetilde{\phi}$ has a continuous extension 
to the boundary and $\widetilde{\phi}=0$ on $b\D$. On a strongly pseudoconvex point 
$p\in b\D$ we have 
\[\phi(p)=\lim_{z\to p}\widetilde{\phi}(z)=0.\] 
Assume that $D$ is a nontrivial analytic disk in $b\D$. Then there exists
a complex line $L\subset \C^2$ such that $\overline{D}=L\cap \Dc$ 
\cite[Lemma 2]{CuckovicSahutoglu09} (see also \cite{FuStraube98}). 

Now we will show that the boundary of $D$ is in the closure of strongly
pseudoconvex points in $b\D$.  Assume that there exists a point $p$ in the
boundary of $D$ that is not in the closure of strongly pseudoconvex points in
$b\D$. Then there exists a nontrivial (affine) analytic disk
$\widetilde{D}\subset b\D$ through $p$ (see, for example, 
\cite[Theorem 1.1]{Freeman74}). Since $\D \subset \C^2$ is convex
we conclude that $\widetilde{D}=\widetilde{L}\cap \Dc$ for a complex line
$\widetilde{L}$. However, $p\in L\cap \widetilde{L}$ and since $\D$ is smooth 
$L$ and $\widetilde{L}$ are perpendicular to the same normal vector at $p$.
Therefore, $L=\widetilde{L}$ and $D=\widetilde{D}$. This contradicts with the
assumption that $p$ is a boundary point of the disk $D$.

The assumption that $\phi\circ f$ is holomorphic for any holomorphic function
$f:\mathbb{D}\to b\D$ implies that the restriction $\phi|_{D}$ is holomorphic
on $D$ for any analytic disk $D\subset b\D$. Then the maximum modulus
principle applied to $\phi|_D$ yields that $\phi=0$ on $D$.  Then continuity of
$\phi$ on $\Dc$ implies that $\phi_1 \phi_2 \cdots \phi_n=\phi=0$ on $b\D$.

To prove iii. $\Rightarrow$ i., again we write  $T=T_{\phi}+K$ where
$\phi=\phi_1\phi_2\cdots\phi_n\in C(\Dc)$ and $K$ is a
compact operator. Then compactness of $T$ follows from the fact that $T_{\phi}$ is 
compact whenever $\phi=0$ on $b\D$.
\end{proof}

\end{document}


\maketitle

We correct an inaccuracy in the statements of Proposition 1, Lemma 1, and 
Corollary 5  in the above paper.  Furthermore, we add a comment and three 
related references.

We start by a definition. A domain $\D$ is said to satisfy condition R if its 
Bergman projection maps $C^{\infty}(\Dc)$ into itself. These domains form 
 a large class. For example, smooth bounded convex domains as 
well as smooth bounded pseudoconvex domains on which the $\dbar$-Neumann 
operator is compact satisfy condition R 
(see \cite[Remark, pg. 82]{BoasStraube91}, \cite{KohnNirenberg65}). However, not 
every smooth bounded pseudoconvex domain satisfies condition R. For example, the 
Diederich-Forn\ae{}ss worm domains do not satisfy condition R \cite{Christ96}.  

The reason for writing this correction is that we use a theorem of Bell, which 
requires the domain to satisfy condition R,  in the proof of the main lemma, 
Lemma \ref{LemDelta} below. In the main results, the condition R is automatic 
since the $\dbar$-Neumann operator is assumed to be compact or the domains are 
convex. However, in Lemma 1, Proposition 1, and Corollary 5 we need to assume 
that the domains satisfy condition R. Now we state the corrected statements.

\begin{lemma} \label{LemDelta}
Let $\D$ be a smooth bounded pseudoconvex domain in $\C^n$ that 
satisfies condition $R, \phi\in C(\Dc),$ and $p\in b\D$ be a strongly 
pseudoconvex point. Then $\lim_{z\to p} \widetilde{\phi}(z)= \phi(p).$  
\end{lemma}

\begin{proposition}
 Let $\D$ be a bounded domain in $\C^n$ and $\psi\in C(\Dc).$
Then $\|T_{\psi}\|_e \leq \|\psi\|_{L^{\infty}(b\D)}.$ Furthermore, 
if $\D$ is smooth bounded pseudoconvex in $\C^n$ that satisfies 
condition $R$ and  the set of strongly pseudoconvex points is dense in 
$b\D$ then $\|T_{\psi}\|_e = \|\psi\|_{L^{\infty}(b\D)}.$
\end{proposition}
\setcounter{corollary}{4}
\begin{corollary}
 Let $\D$ be a smooth bounded pseudoconvex domain in $\C^n$ that
satisfies condition $R$ and the set of strongly pseudoconvex points be 
dense in $b\D$. Assume that $\phi\in C(\Dc)$ and $\lim_{z\to 
p}\widetilde{\phi}(z)=0$ for all strongly pseudoconvex points $p\in b\D$. Then 
$T_{\phi}$ is compact.
\end{corollary}